\documentclass[a4paper, 12pt]{article}
\newcommand{\cqfd}{\hfill\rule{2mm}{2mm}}
\def\as={\,\stackrel{a.s.}{=}\,}
\newdimen\AAdi
\newbox\AAbo

\def\AAk#1#2{\setbox\AAbo=\hbox{#2}\AAdi=\wd\AAbo\kern#1\AAdi{}}
\def\AAr#1#2#3{\setbox\AAbo=\hbox{#2}\AAdi=\ht\AAbo\raise#1\AAdi\hbox{#3}}
\font\tensym=msbm10 \font\sevensym=msbm7 \font\fivesym=msbm5
\newfam\symfam
\textfont\symfam=\tensym \scriptfont\symfam=\sevensym
\scriptscriptfont\symfam=\fivesym

\usepackage[urlcolor=red, 
           linkcolor=blue, 
           colorlinks=red]{hyperref}

\usepackage{amssymb}
\usepackage{eucal}
\usepackage{graphicx}
\usepackage{color}
\usepackage{array,multirow,makecell}
\setcellgapes{1pt}
\makegapedcells
\newcolumntype{R}[1]{>{\raggedleft\arraybackslash }b{#1}}
\newcolumntype{L}[1]{>{\raggedright\arraybackslash }b{#1}}
\newcolumntype{C}[1]{>{\centering\arraybackslash }b{#1}}

\def\1{{1\!\!1}}

\def\as={\,\stackrel{a.s.}{=}\,}

\usepackage[french,english]{babel}
\usepackage[latin1]{inputenc}
\usepackage[T1]{fontenc}
\usepackage{pstricks,pst-node,pst-text,pst-3d}
\usepackage{amsmath}
\usepackage{latexsym}
\usepackage{amssymb}
\usepackage{wasysym}
\usepackage{pifont}
\usepackage{verbatim}
\usepackage{graphicx,dcolumn,babel}
\usepackage{subfigure}
\usepackage{bbm}
\usepackage{amsmath}
\usepackage{amsfonts}
\usepackage[T1]{fontenc}
\newtheorem{theorem}{Theorem}[section]
\newtheorem{lemma}[theorem]{Lemma}
\newtheorem{proposition}[theorem]{Proposition}

\newtheorem{definition}[theorem]{Definition\rm}

\begin{document}
\title{Polynomials shrinkage estimators of a multivariate normal mean}
\author {Abdelkader Benkhaled \footnote{Department of Biology, University of
Mascara, Laboratory of Geomatics, Ecology and Environment (LGEO2E),
Mascara University, Mascara, Algeria. E-mail: benkhaled08@yahoo.fr
},\,{ Mekki Terbeche \footnote{Department of Mathematics,
University of Sciences and Technology, Mohamed Boudiaf, Laboratory
of Analysis and Application of Radiation (LAAR), Oran, USTO-MB,
Algeria. E-mail:mekki.terbeche@gmail.com }} \,and\,{ Abdenour Hamdaoui \footnote{Department of Mathematics, University of Sciences and Technology,
Mohamed Boudiaf, Oran, Laboratory of Statistics and Random
Modelisations of Tlemcen University (LSMA), Algeria. E-mail:
abdenour.hamdaoui@yahoo.fr , abdenour.hamdaoui@unv-usto.dz}
\footnote{Corresponding author: abdenour.hamdaoui@yahoo.fr }}}

\date{}
\maketitle
\begin{abstract}
\noindent In this work, the estimation of the multivariate normal mean by different classes of shrinkage estimators is investigated. The risk associated with the balanced loss function is used to compare two estimators. We start by considering estimators that generalize the James-Stein estimator and show that these estimators dominate the maximum likelihood estimator (MLE), therefore are minimax, when the shrinkage function satisfies some conditions. Then, we treat estimators of polynomial form and prove the increase of the degree of the polynomial allows us to build a better estimator from the one previously constructed.
\end{abstract}

{\small
 }
\section{Introduction}
\label{intro}The multivariate normal distribution has served as a central distribution in much of multivariate analysis. The statistical goal is to estimate the mean parameter which is of interest to many users in almost all fields. The performance of MLE method is not satisfactory, when the dimension of the parameter space is large. The drawbacks of using this method have been shown by Stein \cite{s56} and James and Stein \cite{js61}. Alternative techniques have been developed to improve the MLE; in this paper we focus our attention on shrinkage estimation method. This latter has become a very important technique for modelling data and provides useful techniques for combining data from various sources.
Recent studies, in the context of shrinkage estimation, include Amin
et al.\cite{ana20}, Yuzba et al. \cite{yas20} and Hamdaoui et al.
\cite{hbt20}. Benkhaled and Hamdaoui \cite{bh19}, have considered two forms of shrinkage estimators of the
mean $\theta$ of a multivariate normal distribution $X\sim N_{p}\left(\theta, \sigma^{2}I_{p}\right)$ where
$\sigma^{2}$ is unknown and estimated by the statistic $S^{2}\sim \sigma^{2}\chi_{n}^{2}$. Estimators that shrink
the components of the usual estimator X to zero and estimators of Lindley-type,
that shrink the components of the usual estimator to the random variable X. The
aim is to ameliorate the results of minimaxity obtained in the published papers of estimators cited above.
Hamdaoui et al. \cite{hbm20}, have treated the minimaxity and limits of risks
ratios of shrinkage estimators of a multivariate normal mean in the
Bayesian case. The authors have considered the model $X\sim
N_{p}\left(\theta, \sigma^{2}I_{p}\right)$ where $\sigma^{2}$ is
unknown and have taken the prior law $\theta \sim
N_{p}\left(\upsilon, \tau^{2}I_{p}\right)$. They constructed a
modified Bayes estimator $\delta_{B}^{\ast}$ and an empirical
modified Bayes estimator $\delta_{EB}^{\ast}$. When $n$ and $p$ are
finite, they showed that the estimators $\delta_{B}^{\ast}$ and
$\delta_{EB}^{\ast}$ are minimax. The authors have also interested
in studying the limits of risks ratios of these estimators, to the
MLE $X$, when $n$ and $p$ tend to infinity. The majority of these
authors have been considered the quadratic loss function for computing the risk.\\
A goodness of fit criterion leads to an estimate which gives good
fit and unbiased estimator, thus there is a need to provide a
framework which combines the goodness of fit and precision of
estimation formally. Zellner \cite{z94} suggested balanced losses
that reply this problem. The reader is referred to Guikai et al.
\cite{gqs14}, Karamikabir et al. \cite{kaa18}. Sanjari Farsipour and
Asgharzadeh \cite{sa04} have considered the model: $X_{1},...,X_{n}$
to be a random sample from $N_{p}\left(\theta, \sigma^{2}\right)$
with $\sigma^{2}$ known and the aim is to estimate the parameter
$\theta$. They studied the admissibility of the estimator of the
form $a\overline{X}+b$ under the balanced loss function. Selahattin
and Issam \cite{si19} introduced and derived the optimal extended
balanced loss function (EBLF) estimators and predictors and discussed their performances.\\
In this work, we deal with the model $X\sim N_{p}\left(\theta,
\sigma^{2}I_{p} \right)$, where the parameter $\sigma^{2}$ is
known. Our aim is to estimate the unknown
parameter $\theta$ by shrinkage estimators deduced from the MLE. The adopted
criterion to compare two estimators is the risk
associated to the balanced loss function. The paper is organized as
follows. In Section 2, we recall some preliminaries that are useful
for our main results. In the first part of Section 3, we establish the minimaxity
of the estimators defined by
$\delta_{a}^{(1)}=\left(1-a/\|X\|^{2}\right)X$,
where $\|X\|=(\sum_{i=1}^{p} X_{i}{^2})^{1/2}$ is the euclidean norm of the vector $X=(X_{1},...,X_{p})$ in $\mathrm{R}^{p}$ and the real constant $a$ may depend on $p.$ In the second part of Section 3, we consider the estimators of polynomial form with the
indeterminate $1/\|X\|^{2}$ and show that if we increase the degree of the
polynomial we can build a better estimator from the one previously constructed.
In Section 4, we conduct a simulation study that shows the performance of the considered
estimators. We end the manuscript by giving an Appendix which contains the proofs of some our main results.
\section{Preliminaries}\label{sec:2}
In this section, we recall the following results that are useful in
the proofs of our main results.\\
If $X$ is a multivariate Gaussian random $N_{p}\left(\theta,
\sigma^{2}I_{p}\right) $\textit{\ }in $\mathrm{R}^{p}$, then
$\frac{\Vert X\Vert^{2}}{\sigma^{2}}\sim \chi_{p}^{2}\left(\lambda
\right)$ where $\chi_{p}^{2}\left(\lambda \right)$ denotes the
non-central chi-square distribution with $p$ degrees of freedom and
non-centrality parameter $\lambda=\frac{\left\Vert \theta \right\Vert^{2}}{2\sigma^{2}}$.\\
The following definition given in formula (1.2) by Arnold \cite{a81} will be used to calculate the expectation of
functions of a non-central chi-square law's variable.
\begin{definition} \label{d 2.1}
Let $U \sim \chi_{p}^{2}\left(\lambda \right)$ be non-central chi-square with $p$
degrees of freedom and non-centrality parameter $\lambda$. The density function of $U$ is given by
\begin{eqnarray*}
f(x)=\sum_{k=0}^{+\infty }\frac{e^{-\frac{\lambda}{2}}(\frac{\lambda}{2})^{k}}{k!}\frac{x^{(p/2)+k-1}e^{-x/2}}
{\Gamma(\frac{p}{2}+k)2^{(p/2)+k}}, \ 0<x<+\infty.
\end{eqnarray*}
The right hand side (RHS) of this equality is none other than the formula
\begin{eqnarray*}
\sum_{k=0}^{+\infty }\frac{e^{-\frac{\lambda}{2}}(\frac{\lambda}{2})^{k}}{k!}\chi_{p+2k}^{2},
\end{eqnarray*}
where $\chi_{p+2k}^{2}$ is the density of the central $\chi^{2}$ distribution with $p+2k$ degrees of freedom.
\end{definition}
To this definition we deduce that if $U \sim
\chi_{p}^{2}\left(\lambda \right),$ then for any function $f:
\mathrm{R}_{+} \longrightarrow \mathrm{R}$,
$\chi_{p}^{2}\left(\lambda \right)$ integrable, we have
\begin{eqnarray} \label{e 2}
E\left[f(U)\right]&=& E_{\chi_{p}^{2}\left(\lambda \right)}\left[f(U)\right] \nonumber \\
&=& \int_{\mathrm{R}_{+}}f(x)\chi_{p}^{2}\left(\lambda \right)dx \nonumber\\
&=&\sum_{k=0}^{+\infty }\left[\int_{\mathrm{R}_{+}}f(x)\chi_{p+2k}^{2} dx\right] e^{-\frac{\lambda}{2}}
\frac{\left(\frac{\lambda}{2}\right)^{k}}{k!} \nonumber\\
&=&\sum_{k=0}^{+\infty}\left[\int_{\mathrm{R}_{+}}f(x)\chi_{p+2k}^{2}dx\right]P\left(\frac{\lambda}{2};dk\right),
\end{eqnarray}
where $P\left( \frac{\lambda }{2};dk\right)$ being the Poisson
distribution of parameter $\frac{\lambda }{2}$ and $\chi_{p+2k}^{2}$
is the central chi-square distribution with $p+2k$ degrees of
freedom.\\
The following Stein's Lemma given in \cite{s81} will be often used
in the next.
\begin{lemma} \label{l 2.2}
Let $X$ be a $N\left( \upsilon, \sigma^{2}\right) $ real random
variable and let $f: \mathrm{R} \longrightarrow \mathrm{R}$ be an
indefinite integral of the Lebesgue measurable function,
$f^{\prime}$ essentially the derivative of $f.$ Suppose also that
$E\left(\vert f^{\prime }\left( X\right) \right\vert) <+\infty ,$ then
\begin{eqnarray*}
E\left[\left( \frac{X-\upsilon }{\sigma}\right) f\left( X\right)
\right]=E\left( f^{\prime }\left( X\right) \right).
\end{eqnarray*}
\end{lemma}
\section{Main results}
In this section, we present the model $X\sim N_{p}\left(\theta,\sigma^{2}I_{p} \right)$ where $\sigma^{2}$ is known. Our aim is to estimate the unknown mean parameter $\theta$ by the shrinkage estimators under the balanced squared error loss function. For the sake of simplicity, we treat only the case when $\sigma^{2}=1$, as long as by a change of variable, any model of type $Y\sim N_{p}\left(\theta_{1} ,\sigma^{2}I_{p} \right)$ can be reduced to the model $Z\sim N_{p}\left(\theta_{2} ,I_{p} \right)$. Namely, we consider the model $X\sim N_{p}\left(\theta ,I_{p} \right)$ and we want to estimate the unknown parameter $\theta$.
\begin{definition} \label{d 3.1}
Suppose that $X$ is a random vector having a multivariate normal distribution $N_{p}\left(\theta ,
I_{p} \right)$ where the parameter $\theta$ is unknown. The balanced squared error loss function is defined as follows:
\begin{eqnarray} \label{e 3}
L_{\omega}(\delta,\theta)=\omega
\|\delta-\delta_{0}\|^{2}+(1-\omega)\|\delta-\theta\|^{2}, \ 0\leq \omega<1,
\end{eqnarray}
where $\delta_{0}$ is the target estimator of $\theta$, $\omega$ is the weight given to the proximity
of $\delta$ to $\delta_{0}$, $1-\omega$ is the relative weight given to the precision of estimation
portion and $\delta$ is a given estimator.
\end{definition}
For more details about this loss see Jafari Jozani et al.
\cite{jlm14}, Zinodiny et al. \cite{zln17} and Karamikabir and
Afsahri \cite{ka20}.\\
We associate to this balanced squared error loss function the risk
function defined by $$R_{\omega}(\delta,\theta)=E(L_{\omega}(\delta,\theta)).$$
In this model, it is clear that the MLE is $X:=\delta_{0}$, its risk function is $(1-\omega)p$.\\
Indeed: we have
\begin{eqnarray*}
R_{\omega}(X,\theta)&=&\omega E(\|X-X\|^{2})+(1-\omega)E(\|X-\theta\|^{2})\\
&=&(1-\omega)E(\|X-\theta\|^{2}).
\end{eqnarray*}
As $X\sim N_{p}\left(\theta ,I_{p} \right)$, then
$X-\theta\sim N_{p}\left(0 ,I_{p} \right)$, therefore
$\|X-\theta\|^{2}\sim \chi^{2}_{p}$.\\
Hence, $E(\|X-\theta\|^{2})=E(\chi^{2}_{p})=p,$
and the desired result follows.\\
It is well known that $\delta_{0}$ is minimax and inadmissible for
$p\geq3$, thus any estimator dominates it is also minimax. We give
the following Lemma, that will be used in our proofs and its proof
is postponed to the Appendix.
\begin{lemma} \label{l 3.1}
Let $U\sim\chi_{p}^{2}(\lambda)$ be non-central chi-square with $p$
degrees of freedom and non-centrality parameter $\lambda$ then,
\begin{itemize}
\item[i)] for any real numbers $s$ and $r$ where $-\frac{p}{2}<s \leq r<0,$ the real function
\begin{eqnarray*}
H_{p,r,s}(\lambda)=\frac{E(U^{r})}{E(U^{s})}=\frac{\int_{R_{+}}x^{r}\chi_{p}^{2}(\lambda;dx)}
{\int_{R_{+}}x^{s}\chi_{p}^{2}(\lambda;dx)}
\end{eqnarray*}
is nondecreasing on $\lambda.$
\item[ii)] Furthermore, if $X\sim N_{p}\left(\theta ,I_{p} \right)$, we get
\begin{eqnarray*}
\sup_{\|\theta\|}\left(\frac{E(\|X\|^{-2r+2})}{E(\|X\|^{-r})}\right)=2^{\frac{-r+2}{2}}
\frac{\Gamma(\frac{p}{2}-r+1)}{\Gamma(\frac{p-r}{2})}.
\end{eqnarray*}
\end{itemize}
\end{lemma}

\subsection{James-Stein estimators and minimaxity}
In 1956, Stein \cite{s56} proved a result that astonished many researchers and was catalyst an enormous and rich literature of substantial importance in statistical theory and practice. He showed that when estimating, under squared error loss, the unknown mean vector $\theta$ of a $p$-dimensional random vector $X$ having a normal distribution with identity covariance matrix, estimators of the form $\delta_{a,b}=\left(1-a/(b+\|X\|^{2})\right)X$ dominate the usual estimator $X$ for $a$ sufficiently small and $b$ sufficiently large when $p \geq 3$. In 1961, James and Stein \cite{js61} sharpened the result and gave an explicit class of dominating estimators, $\delta_{a}=\left(1-a/\|X\|^{2})\right)X$ for $0<a<2(p-2)$, and also showed that the choice on $a=p-2$ (the James-Stein estimator) is uniformly best. In this section we show the sufficient condition for which the estimator $\delta_{a}$ dominates the usual estimator $X$ under the balanced loss function $L_{\omega}$ defined in (\ref{e 3}) and we determined the optimal value for $a$ (corresponding to the James-Stein estimator) that minimizes the risk function $R_{\omega}(\delta_{a},\theta)$. 

Consider the estimator
\begin{eqnarray} \label{e 3.11}
\delta_{a}^{(1)}=\left(1-\frac{a}{\|X\|^{2}}\right)X
=X-\frac{a}{\|X\|^{2}}X,
\end{eqnarray}
where the real constant $a$ may depend on $p.$
\begin{proposition} \label{p 3.11}
Under the balanced loss function $L_{\omega}$, the risk function of
the estimator $\delta_{a}^{(1)}$ given in (\ref{e 3.11}) is
\begin{eqnarray*}
R_{\omega}(\delta_{a}^{(1)},\theta)&=&(1-\omega)\left[p-2a(p-2)E\left(\frac{1}{\|X\|^{2}}\right)\right]
+a^{2}E\left(\frac{1}{\|X\|^{2}}\right).
\end{eqnarray*}
\end{proposition}
\paragraph{Proof} Using the risk function associated to the balanced loss function $L_{\omega}$ defined in 
in (\ref{e 3}) and the formula of the estimator $\delta_{a}^{(1)}$ given in (\ref{e 3.11}), we obtain
\begin{eqnarray*}
R_{\omega}(\delta_{a}^{(1)},\theta)&=&\omega
E\left(\left\|-\frac{a}{\|X\|^{2}}X\right\|^{2}\right)
+(1-\omega)E\left(\left\|X-\theta-\frac{a}{\|X\|^{2}}X\right\|^{2}\right)\\
&=&a^{2}E\left(\frac{1}{\|X\|^{2}}\right)+(1-\omega)p-2a(1-\omega)E\left(\left\langle
X-\theta,\frac{1}{\|X\|^{2}}X\right\rangle\right).
\end{eqnarray*}
As,
\begin{eqnarray*}
E\left(\left\langle
X-\theta,\frac{1}{\|X\|^{2}}X\right\rangle\right)&=&\sum^{p}_{i=1}
E\left[(X_{i}-\theta_{i})\frac{1}{\|X\|^{2}}X_{i}\right]
\end{eqnarray*}
Using Lemma \ref{l 2.2}, we get
\begin{eqnarray*}
E\left(\left\langle X-\theta,\frac{1}{\|X\|^{2}}X\right\rangle\right)&=&\sum^{p}_{i=1}
E\left(\frac{\partial}{\partial X_{i}}\frac{1}{\|X\|^{2}}X_{i}\right)\\
&=&\sum^{p}_{i=1}E\left(\frac{1}{\|X\|^{2}}-\frac{2 X_{i}^{2}}{\|X\|^{4}}\right)\\
&=&(p-2)E\left(\frac{1}{\|X\|^{2}}\right).
\end{eqnarray*}
Then
\begin{eqnarray*}
R_{\omega}(\delta_{a}^{(1)},\theta)&=&a^{2}E\left(\frac{1}{\|X\|^{2}}\right)+(1-\omega)p
-2a(1-\omega)E\left(\left\langle X-\theta,\frac{1}{\|X\|^{2}}X\right\rangle\right)\\
&=&a^{2}E\left(\frac{1}{\|X\|^{2}}\right)+(1-\omega)p
-2a(1-\omega)(p-2)E\left(\frac{1}{\|X\|^{2}}\right)\\
&=&(1-\omega)\left[p-2a(p-2)E\left(\frac{1}{\|X\|^{2}}\right)\right]+a^{2}E\left(\frac{1}{\|X\|^{2}}\right).
\end{eqnarray*}\cqfd\\
Using the convexity on $a$ of the function
$R_{\omega}(\delta_{a}^{(1)},\theta)$, the optimal value for $a$
that minimizes the risk function
$R_{\omega}(\delta_{a}^{(1)},\theta)$, is
\begin{eqnarray} \label{e 3.111}
\widehat a=(1-\omega)(p-2).
\end{eqnarray}
For $a=\widehat a$, we obtain the James-Stein estimator
\begin{eqnarray} \label{e 3.12}
\delta_{JS}=\delta_{\widehat a,2}=\left(1-\frac{\widehat
a}{\|X\|^{2}}\right)X.
\end{eqnarray}
From Proposition \ref {p 3.11}, the risk function of $\delta_{JS}$ is
\begin{eqnarray}\label{e 3.13}
R_{\omega}(\delta_{JS},\theta)=(1-\omega)p-(p-2)^{2}(1-\omega)^{2}
E\left(\frac{1}{p-2+2K}\right),
\end{eqnarray}
where $K\sim P\left(\frac{\|\theta\|^{2}}{2}\right)$.\\
From the formula (\ref{e 3.13}), we note that
$R_{\omega}(\delta_{JS},\theta)\leq(1-\omega)p=R_{\omega}(X,\theta),$
then $\delta_{JS}$ dominates the MLE $X,$ therefore it is also
minimax.
\subsection{Polynomials shrinkage estimators}
Since the estimator $\delta_{a}=X-a\frac{1}{\|X\|^{2}}X$ dominates the MLE $X$ for certain values of $a$, we think to add the term $b(\frac{1}{\|X\|^{2}})^{2} X$ to the James-Stein estimator $\delta_{JS}$ to obtain an estimator that outperforms $\delta_{JS}$, then we construct the classes of shrinkage estimators which dominate the James-Stein estimator $\delta_{JS}$. Our main idea is to add each time a term of the form $\gamma (1/\|X\|^{2})^{m} X$ where $\gamma$ is a real constant may depend on $p$ and the parameter $m$ is integer, and we construct the estimators which dominate the estimators of the class defined previously. Thus in this section we deal with the shrinkage estimators of polynomial form with the indeterminate $1/\|X\|^{2}$.

Let the estimator
\begin{eqnarray} \label{e 3.21}
\delta_{b}^{(2)}&=&\delta_{JS}+b\left(\frac{1}{\|X\|^{2}}\right)^{2}X \nonumber\\
&=&X-(1-\omega)(p-2)\frac{1}{\|X\|^{2}}X+b\left(\frac{1}{\|X\|^{2}}\right)^{2}X,
\end{eqnarray}
where the real constant $b$ may depend on $p.$
\begin{proposition} \label{p 3.21}
Under the balanced loss function $L_{\omega}$, the risk function of
the estimator $\delta_{b}^{(2)}$ given in (\ref{e 3.21}) is
\begin{eqnarray*}
R_{\omega}(\delta_{b}^{(2)},\theta)&=&R_{\omega}(\delta_{JS},\theta)-4b(1-\omega)
E\left(\frac{1}{\|X\|^{4}}\right)+b^{2}E\left(\frac{1}{\|X\|^{6}}\right).
\end{eqnarray*}
\end{proposition}
\paragraph{Proof} Using the risk function associated to the balanced loss function $L_{\omega}$ defined in 
in (\ref{e 3}) and the formula of the estimator $\delta_{b}^{(2)}$ given in (\ref{e 3.21}), we get
\begin{eqnarray*}
\hspace{-1cm} R_{\omega}(\delta_{b}^{(2)},\theta)&=&\omega
E\left(\left\|\delta_{JS}+b\frac{1}{\|X\|^{2}}X-X\right\|^{2}\right)+(1-\omega)E\left(\left\|\delta_{JS}+
b\frac{1}{\|X\|^{2}}X-\theta\right\|^{2}\right)\\
&=&\omega
E\left(\left\|\delta_{JS}-X\right\|^{2}+b^{2}\frac{1}{(\|X\|^{2})^{3}}+2\left\langle
\delta_{JS}-X,b\frac{1}{(\|X\|^{2})^{2}}X\right\rangle\right)\\
&+&(1-\omega)E\left(\left\|\delta_{JS}-\theta\right\|^{2}+b^{2}\frac{1}{(\|X\|^{2})^{3}}+2\left\langle
\delta_{JS}-\theta,b\frac{1}{(\|X\|^{2})^{2}}X\right\rangle\right)
\end{eqnarray*}
\begin{eqnarray*}
&=&R_{\omega}(\delta_{JS},\theta)+b^{2}E\left(\frac{1}{(\|X\|^{2})^{3}}\right)
-2b\omega(1-\omega)(p-2)E\left(\frac{1}{(\|X\|^{2})^{2}}\right)\\
&+&2(1-\omega)E\left(\left\langle
X-\theta-(1-\omega)(p-2)\frac{1}{\|X\|^{2}}X,
b\frac{1}{(\|X\|^{2})^{2}}X\right\rangle\right)\\
&=&R_{\omega}(\delta_{JS},\theta)+b^{2}E\left(\frac{1}{\|X\|^{6}}\right)
-2b\omega(1-\omega)(p-2)E\left(\frac{1}{\|X\|^{4}}\right)\\
&+&2b(1-\omega)\sum^{p}_{i=1}E\left((X_{i}-\theta_{i})\frac{X_{i}}{\|X\|^{4}}\right)
-2b(1-\omega)^{2}(p-2)E\left(\frac{1}{\|X\|^{4}}\right).
\end{eqnarray*}
Using Lemma \ref{l 2.2}, we get
\begin{eqnarray*}
\sum^{p}_{i=1}E\left((X_{i}-\theta_{i})\frac{X_{i}}{\|X\|^{4}}\right)&=&\sum^{p}_{i=1}E\left(\frac{\partial}{\partial
X_{i}}\frac{X_{i}}{\|X\|^{4}}\right)\\
&=&\sum^{p}_{i=1}E\left(\frac{1}{\|X\|^{4}}-4\frac{X_{i}^{2}}{\|X\|^{6}}\right)\\
&=&(p-4)E\left(\frac{1}{\|X\|^{4}}\right).
\end{eqnarray*}
Then,
\begin{eqnarray*}
\hspace{-1cm}R_{\omega}(\delta_{b}^{(2)},\theta)&=&R_{\omega}(\delta_{JS},\theta)+b^{2}
E\left(\frac{1}{\|X\|^{6}}\right)-2b\omega(1-\omega)(p-2)E\left(\frac{1}{\|X\|^{4}}\right)\\
&+&2b(1-\omega)(p-4)E\left(\frac{1}{\|X\|^{4}}\right)-2b(1-\omega)^{2}(p-2)E\left(\frac{1}{\|X\|^{4}}\right)\\
&=&R_{\omega}(\delta_{JS},\theta)-4b(1-\omega)
E\left(\frac{1}{\|X\|^{4}}\right)+b^{2}E\left(\frac{1}{\|X\|^{6}}\right).
\end{eqnarray*}\cqfd\\
\begin{theorem} \label{t 3.22}
Under the balanced loss function $L_{\omega}$, the estimator
$\delta_{b}^{(2)}$ with $p>6$ and
$$b=2(1-\omega)(p-6),$$ dominates
the James-Stein estimator $\delta_{JS}.$
\end{theorem}
\paragraph{Proof} Using the Proposition \ref{p 3.21}, we have
\begin{eqnarray*}
R_{\omega}(\delta_{b}^{(2)},\theta)&=&R_{\omega}(\delta_{JS},\theta)-4b(1-\omega)
E\left(\frac{1}{\|X\|^{4}}\right)+b^{2}\frac{E\left(\frac{1}{\|X\|^{6}}\right)}
{E\left(\frac{1}{\|X\|^{4}}\right)}E\left(\frac{1}{\|X\|^{4}}\right).
\end{eqnarray*}
From ii) of Lemma \ref{l 3.1}, we obtain
\begin{eqnarray*}
\frac{E\left(\frac{1}{\|X\|^{6}}\right)}
{E\left(\frac{1}{\|X\|^{4}}\right)}=\frac{E\left(\|X\|^{-6}\right)}
{E\left(\|X\|^{-4}\right)}\leq
2^{\frac{-4+2}{2}}\frac{\Gamma(\frac{p}{2}-4+1)}{\Gamma(\frac{p-4}{2})}=\frac{1}{p-6}.
\end{eqnarray*}
Then,
\begin{eqnarray} \label{e 3.22}
R_{\omega}(\delta_{b}^{(2)},\theta)&\leq&R_{\omega}(\delta_{JS},\theta)-
4b(1-\omega)E\left(\frac{1}{\|X\|^{4}}\right)+b^{2}\frac{1}{(p-6)}
E\left(\frac{1}{\|X\|^{4}}\right).
\end{eqnarray}
The optimal value for $b$ that minimizes the right hand side of the last
inequality, is
\begin{eqnarray} \label{e 3.221}
\widehat b=2(1-\omega)(p-6).
\end{eqnarray}
If we replace $b$ by $\widehat b$ in the inequality (\ref{e 3.22}), we get
\begin{eqnarray*}
R_{\omega}(\delta_{\widehat b}^{(2)},\theta)&\leq&
R_{\omega}(\delta_{JS},\theta)-4(1-\omega)^{2}(p-6)E\left(\frac{1}{\|X\|^{4}}\right)\\
&\leq& R_{\omega}(\delta_{JS},\theta).
\end{eqnarray*}\cqfd\\
Now, we consider the estimator
\begin{eqnarray} \label{e 3.31}
\hspace{-1cm}\delta_{c}^{(3)}&=&\delta_{\widehat b}^{(2)}+c\left(\frac{1}{\|X\|^{2}}\right)^{3}X \nonumber\\
&=&X-\widehat a \frac{1}{\|X\|^{2}}X+ \widehat b
\left(\frac{1}{\|X\|^{2}}\right)^{2}X+c\left(\frac{1}{\|X\|^{2}}\right)^{3}X,
\end{eqnarray}
where the constants $\widehat a$ and $\widehat b$ are given
respectively in (\ref{e 3.111}) and (\ref{e 3.221}) and the real
parameter $c$ may depend on $p.$
\begin{proposition} \label{p 3.31}
Under the balanced loss function $L_{\omega}$, the risk function of
the estimator $\delta_{c}^{(3)}$ given in (\ref{e 3.31}) is
\begin{eqnarray*}
R_{\omega}(\delta_{c}^{(3)},\theta)&=&R_{\omega}(\delta_{\widehat
b}^{(2)},\theta)+c^{2}E\left(\frac{1}{\|X\|^{10}}\right)+4c(1-\omega)(p-6)E\left(\frac{1}{\|X\|^{8}}\right)\\
&-&8c(1-\omega)E\left(\frac{1}{\|X\|^{6}}\right).
\end{eqnarray*}
\end{proposition}
\paragraph{Proof} Using the risk function associated to the balanced loss function $L_{\omega}$ defined in 
in (\ref{e 3}) and the formula of the estimator $\delta_{c}^{(3)}$ given in (\ref{e 3.31}), we obtain
\begin{eqnarray*}
R_{\omega}(\delta_{c}^{(3)},\theta)&=&\omega E\left(\left\|\delta_{\widehat
b}^{(2)}+c\left(\frac{1}{\|X\|^{2}}\right)^{3}X-X\right\|^{2}\right)\\
&+&(1-\omega)E\left(\left\|\delta_{\widehat
b}^{(2)}+c\left(\frac{1}{\|X\|^{2}}\right)^{3}X-\theta\right\|^{2}\right)\\
&=&\omega E\left(\left\|\delta_{\widehat b}^{(2)}-X\right\|^{2}+c^{2}\frac{1}
{(\|X\|^{2})^{5}}+2\left\langle\delta_{\widehat b}^{(2)}-X,c\frac{1}{(\|X\|^{2})^{3}}X\right\rangle\right)\\
&+&(1-\omega)E\left(\left\|\delta_{\widehat b}^{(2)}-\theta\right\|^{2}+c^{2}\frac{1}{(\|X\|^{2})^{5}}
+2\left\langle\delta_{\widehat b}^{(2)}-\theta,c\frac{1}{(\|X\|^{2})^{3}}X\right\rangle\right)\\
&=&R_{\omega}(\delta_{\widehat b}^{(2)},\theta)+c^{2}E\left(\frac{1}{(\|X\|^{2})^{5}}\right)\\
&+&2\omega E\left\langle -\widehat a\frac{1}{\|X\|^{2}}X+\widehat b\left(\frac{1}{\|X\|^{2}}\right)^{2}X,
c\left(\frac{1}{\|X\|^{2}}\right)^{3}X\right\rangle\\
&+&2(1-\omega)E \left\langle X-\theta- \widehat a\frac{1}{\|X\|^{2}}X+\widehat b\left(\frac{1}
{\|X\|^{2}}\right)^{2}X,c\left(\frac{1}{\|X\|^{2}}\right)^{3}X\right\rangle\\
&=&R_{\omega}(\delta_{\widehat b}^{(2)},\theta)+c^{2}E\left(\frac{1}{\|X\|^{10}}\right)-2c\widehat
aE\left(\frac{1}{\|X\|^{6}}\right)\\
&+&2c\widehat b
E\left(\frac{1}{\|X\|^{8}}\right)+2c(1-\omega)\sum^{p}_{i=1}E\left((X_{i}-\theta_{i})\frac{X_{i}}{\|X\|^{6}}\right).
\end{eqnarray*}
Using Lemma \ref{l 2.2}, we get
\begin{eqnarray*}
\sum^{p}_{i=1}E\left((X_{i}-\theta_{i})\frac{X_{i}}{\|X\|^{6}}\right)&=&
\sum^{p}_{i=1}E\left(\frac{\partial}{\partial
X_{i}}\frac{X_{i}}{\|X\|^{6}}\right)\\
&=&(p-6)E\left(\frac{1}{\|X\|^{6}}\right).
\end{eqnarray*}
Then
\begin{eqnarray*}
\hspace{-1cm}R_{\omega}(\delta_{c}^{(3)},\theta)&=&R_{\omega}(\delta_{\widehat
b}^{(2)},\theta)+c^{2}E\left(\frac{1}{\|X\|^{10}}\right)-2c(1-\omega)(p-2)E\left(\frac{1}{\|X\|^{6}}\right)\\
&+&4c(1-\omega)(p-6)E\left(\frac{1}{\|X\|^{8}}\right)+2c(1-\omega)(p-6)E\left(\frac{1}{\|X\|^{6}}\right)\\
&=&R_{\omega}(\delta_{\widehat
b}^{(2)},\theta)+c^{2}E\left(\frac{1}{\|X\|^{10}}\right)+4c(1-\omega)(p-6)E\left(\frac{1}{\|X\|^{8}}\right)\\
&-&8c(1-\omega)E\left(\frac{1}{\|X\|^{6}}\right).
\end{eqnarray*}\cqfd\\
\begin{theorem} \label{t 3.32}
Under the balanced loss function $L_{\omega}$, the estimator
$\delta_{c}^{(3)}$ with $p>10$ and
$$c=2(1-\omega)(p-10)^{2},$$ dominates
the estimator $\delta_{\widehat b}^{(2)}.$
\end{theorem}
\paragraph{Proof} Using the last Proposition, we have
\begin{eqnarray*}
\hspace{-1cm}R_{\omega}(\delta_{c}^{(3)},\theta)&=&
R_{\omega}(\delta_{\widehat
b}^{(2)},\theta)+c^{2}\frac{E\left(\frac{1}{\|X\|^{10}}\right)}
{E\left(\frac{1}{\|X\|^{6}}\right)}E\left(\frac{1}{\|X\|^{6}}\right)\\
&+&4c(1-\omega)(p-6)\frac{E\left(\frac{1}{\|X\|^{8}}\right)}
{E\left(\frac{1}{\|X\|^{6}}\right)}E\left(\frac{1}{\|X\|^{6}}\right)\\
&-&8c(1-\omega)E\left(\frac{1}{\|X\|^{6}}\right).
\end{eqnarray*}
From ii) of Lemma \ref{l 3.1}, we obtain
\begin{eqnarray*}
\frac{E\left(\frac{1}{\|X\|^{10}}\right)}
{E\left(\frac{1}{\|X\|^{6}}\right)}=\frac{E\left(\|X\|^{-10}\right)}
{E\left(\|X\|^{-6}\right)}\leq
2^{\frac{-6+2}{2}}\frac{\Gamma(\frac{p}{2}-6+1)}{\Gamma(\frac{p-6}{2})}=\frac{1}{(p-8)(p-10)},
\end{eqnarray*}
and from i) of Lemma \ref{l 3.1}, we get
\begin{eqnarray*}
\frac{E\left(\frac{1}{\|X\|^{8}}\right)}
{E\left(\frac{1}{\|X\|^{6}}\right)}&=&\frac{E\left((\chi_{p}^{2}(\lambda))^{-4}\right)}{E\left((\chi_{p}^{2}(\lambda))^{-3}\right)}\\
&\leq&\frac{E\left((\chi_{p}^{2})^{-4}\right)}{E\left((\chi_{p}^{2})^{-3}\right)}
=\frac{2^{-4}\frac{\Gamma(\frac{p}{2}-4)}{\Gamma(\frac{p}{2})}}{2^{-3}\frac{\Gamma(\frac{p}{2}-3)}{\Gamma(\frac{p}{2})}}
=\frac{1}{p-8}.
\end{eqnarray*}
where $\lambda=\frac{\|\theta\|^{2}}{2}$ and $\chi_{p}^{2}$
is the central chi-square distribution with $p$ degrees of
freedom. Then,
\begin{eqnarray}\label{e 3.32}
R_{\omega}(\delta_{c}^{(3)},\theta)&\leq&R_{\omega}(\delta_{\widehat
b}^{(2)},\theta)
+c^{2}\frac{1}{(p-8)(p-10)}E\left(\frac{1}{\|X\|^{6}}\right)\nonumber\\
&+&4c(1-\omega)(p-6)\frac{1}{p-8}E\left(\frac{1}{\|X\|^{6}}\right)-8c(1-\omega)E\left(\frac{1}{\|X\|^{6}}\right)\nonumber\\
&=& R_{\omega}(\delta_{\widehat b}^{(2)},\theta)+c^{2}\frac{1}{(p-8)(p-10)}E\left(\frac{1}{\|X\|^{6}}\right)\nonumber\\
&-&4c(1-\omega)\frac{p-10}{p-8}E\left(\frac{1}{\|X\|^{6}}\right).
\end{eqnarray}
The optimal value for $c$ that minimizes the right hand side of the
inequality (\ref{e 3.32}), is
\begin{eqnarray} \label{e 3.331}
\widehat c=2(1-\omega)(p-10)^{2}.
\end{eqnarray}
If we replace $c$ by $\widehat c$ in the inequality (\ref{e 3.32}), we get
\begin{eqnarray*}
R_{\omega}(\delta_{c}^{(3)},\theta)&\leq&R_{\omega}(\delta_{\widehat
b}^{(2)},\theta)-4\frac{(1-\omega)^{2}(p-10)^{3}}{p-8}E\left(\frac{1}{\|X\|^{6}}\right)\\
&\leq& R_{\omega}(\delta_{\widehat b}^{(2)},\theta).
\end{eqnarray*}\cqfd\\
Now, we consider the estimator
\begin{eqnarray} \label{e 3.41}
\delta_{d}^{(4)}&=&\delta_{\widehat c}^{(2)}+d\left(\frac{1}{\|X\|^{2}}\right)^{4}X \nonumber\\
&=&X-\widehat a \frac{1}{\|X\|^{2}}X+ \widehat b
\left(\frac{1}{\|X\|^{2}}\right)^{2}X+\widehat
c\left(\frac{1}{\|X\|^{2}}\right)^{3}X+d\left(\frac{1}{\|X\|^{2}}\right)^{4}X,
\end{eqnarray}
where the constants $\widehat a$ and $\widehat b$ and $\widehat c$
are given respectively in (\ref{e 3.111}), (\ref{e 3.221}) and
(\ref{e 3.331}) and the real parameter $d$ may depend on $p$.
Using the same technique used in the proofs of Proposition (\ref{p 3.31})
and theorem (\ref{t 3.32}), we obtain the following results.
\begin{proposition} \label{p 3.41}
Under the balanced loss function $L_{\omega}$, the risk function of
the estimator $\delta_{d}^{(4)}$ given in (\ref{e 3.41}) is
\begin{eqnarray*}
R_{\omega}(\delta_{d}^{(4)},\theta)&=&R_{\omega}(\delta_{\widehat c}^{(3)},\theta)+d^{2}E\left(\frac{1}
{\|X\|^{14}}\right)+4d(1-\omega)(p-10)^{2}E\left(\frac{1}{\|X\|^{12}}\right)\\
&+&4d(1-\omega)(p-6)E\left(\frac{1}{\|X\|^{10}}\right)-12d(1-\omega)E\left(\frac{1}{\|X\|^{8}}\right),
\end{eqnarray*}
\end{proposition}
\begin{theorem} \label{t 3.42}
Under the balanced loss function $L_{\omega}$, the estimator
$\delta_{d}^{(4)}$ with $p>14$ and
$$ d=2(1-\omega)(p^{2}-28p+188)(p-14) $$
dominates the estimator $\delta_{\widehat c}^{(3)}$.
\end{theorem}
\section{Simulation results}
We recall the form of the James-Stein estimator $\delta_{JS}$ given
in (\ref{e 3.12}). Its risk function associated to the balanced squared error loss function 
$L_{\omega}$ is given by the formula (\ref{e 3.13}). \\
We also recall the estimators $\delta_{b}^{(2)},$ $\delta_{c}^{(3)},$ and $\delta_{d}^{(4)},$ given respectively in (\ref{e 3.21}), (\ref{e 3.31}) and (\ref{e 3.41})
with $b=(1-\omega)(p-6),$ $c=(1-\omega)(p-10)^{2}$ and $d=2(1-\omega)(p^{2}-28p+188)(p-14)$. Their risk functions associated to the balanced squared error loss function $L_{\omega}$ are obtained by replacing the constants $b, c$ and $d$ in the Propositions (\ref{p 3.21}), (\ref{p 3.31}) and (\ref{p 3.41}) respectively.    

In this section, taking the values of the constants $b,$ $c$ and $d$ given above. In the first part of this section, we present the graphs of the risks ratios of
the estimators $\delta_{JS},$ $\delta_{b}^{(2)}$ and $\delta_{c}^{(3)}$, to the MLE $X$ denoted respectively:
$\frac{R_{\omega}(\delta_{JS},\theta)}{R_{\omega}(X,\theta)},$
$\frac{R_{\omega}(\delta_{b}^{(2)},\theta)}{R_{\omega}(X,\theta)}$ and
$\frac{R_{\omega}(\delta_{c}^{(3)},\theta)}{R_{\omega}(X,\theta)}$ as
function of $\lambda=\left\Vert \theta
\right\Vert^{2},$ for various values of $p$ and
$\omega.$ In the second part of this section, we present two groups of tables. The first group containing the values of
risks ratios
$\frac{R_{\omega}(\delta_{JS},\theta)}{R_{\omega}(X,\theta)},$
$\frac{R_{\omega}(\delta_{b}^{(2)},\theta)}{R_{\omega}(X,\theta)}$ and $\frac{R_{\omega}(\delta_{c}^{(3)},\theta)}{R_{\omega}(X,\theta)}$ as a function of variable $\lambda=\left\Vert \theta
\right\Vert^{2},$ for various values of $p$ and
$\omega.$ In the second group we give the values of
risks ratios $\frac{R_{\omega}(\delta_{c}^{(3)},\theta)}{R_{\omega}(X,\theta)}$ and $\frac{R_{\omega}(\delta_{d}^{(4)},\theta)}{R_{\omega}(X,\theta)}$ as a function of variable $\lambda=\left\Vert \theta
\right\Vert^{2},$ for various values of $p$ and
$\omega.$ 
\begin{figure}[htb]
 \centering
 \includegraphics[height=7.1cm,width=10.0cm]{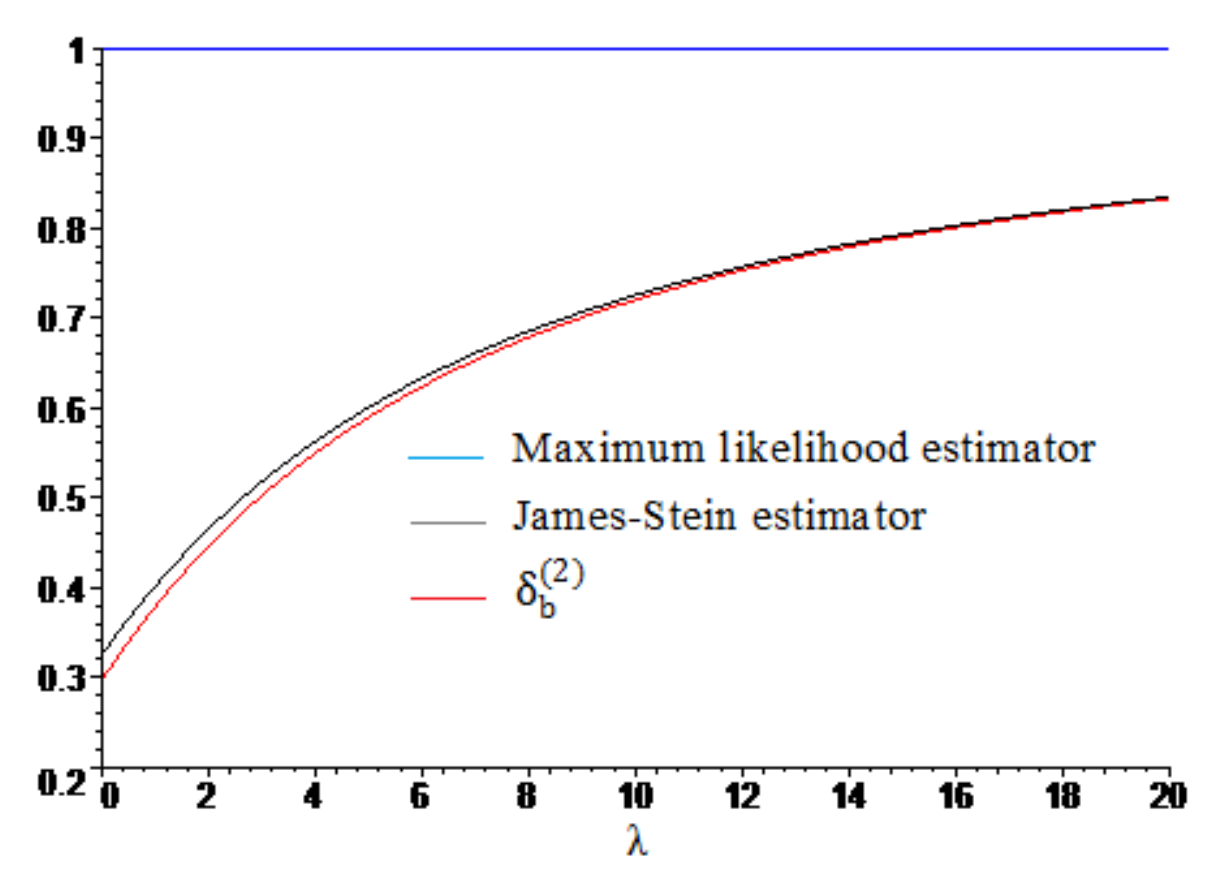}
 \caption{\rm Graph of risks ratios
$\frac{R_{\omega}(\delta_{JS},\theta)}{R_{\omega}(X,\theta)}$ and
$\frac{R_{\omega}(\delta_{b}^{(2)},\theta)}{R_{\omega}(X,\theta)}$ as function of
$\lambda$ for $p=8$ and $\omega=0.1$ } \label{figure:f1}
\end{figure}

\newpage

\begin{figure}[htb]
 \centering
 \includegraphics[height=7.1cm,width=10.0cm]{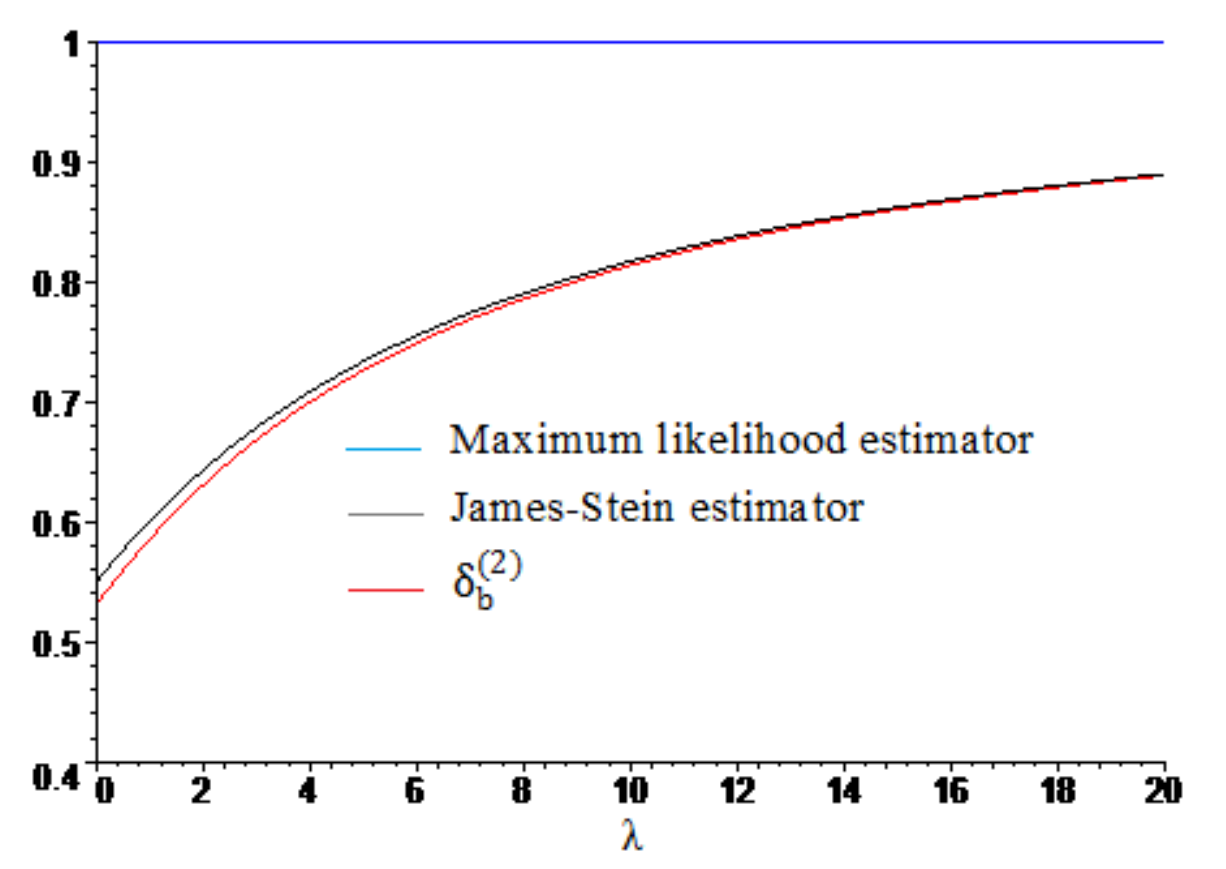}
 \caption{\rm Graph of risks ratios
$\frac{R_{\omega}(\delta_{JS},\theta)}{R_{\omega}(X,\theta)}$ and
$\frac{R_{\omega}(\delta_{b}^{(2)},\theta)}{R_{\omega}(X,\theta)}$ as function of
$\lambda$ for $p=8$ and $\omega=0.4$} \label{figure:f2}
%
 \centering
 \includegraphics[height=7.1cm,width=10.0cm]{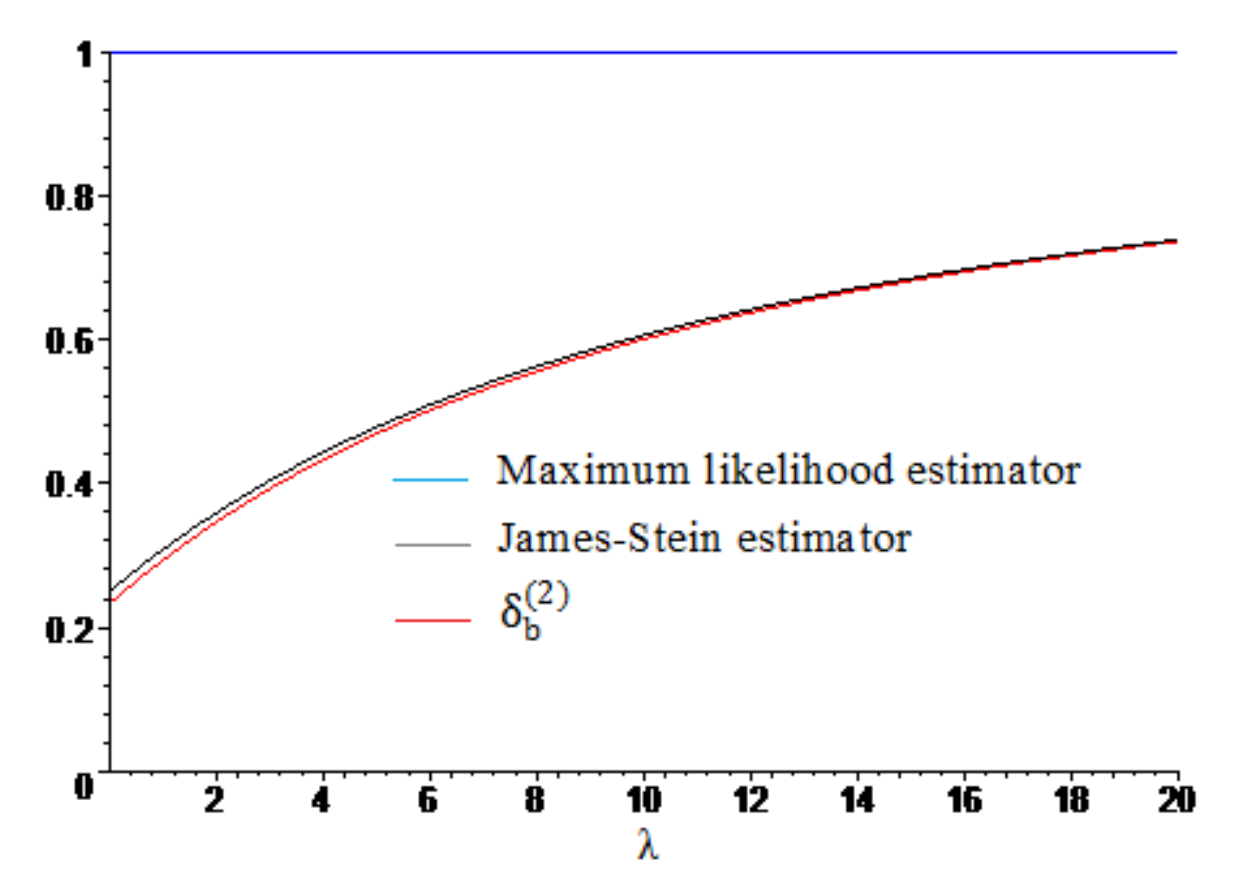}
 \caption{\rm Graph of risks ratios
$\frac{R_{\omega}(\delta_{JS},\theta)}{R_{\omega}(X,\theta)}$ and
$\frac{R_{\omega}(\delta_{b}^{(2)},\theta)}{R_{\omega}(X,\theta)}$ as function of
$\lambda$ for $p=12$ and $\omega=0.1$ } \label{figure:f3}
\end{figure}

\newpage

\begin{figure}[htb]
 \centering
 \includegraphics[height=7.1cm,width=10.0cm]{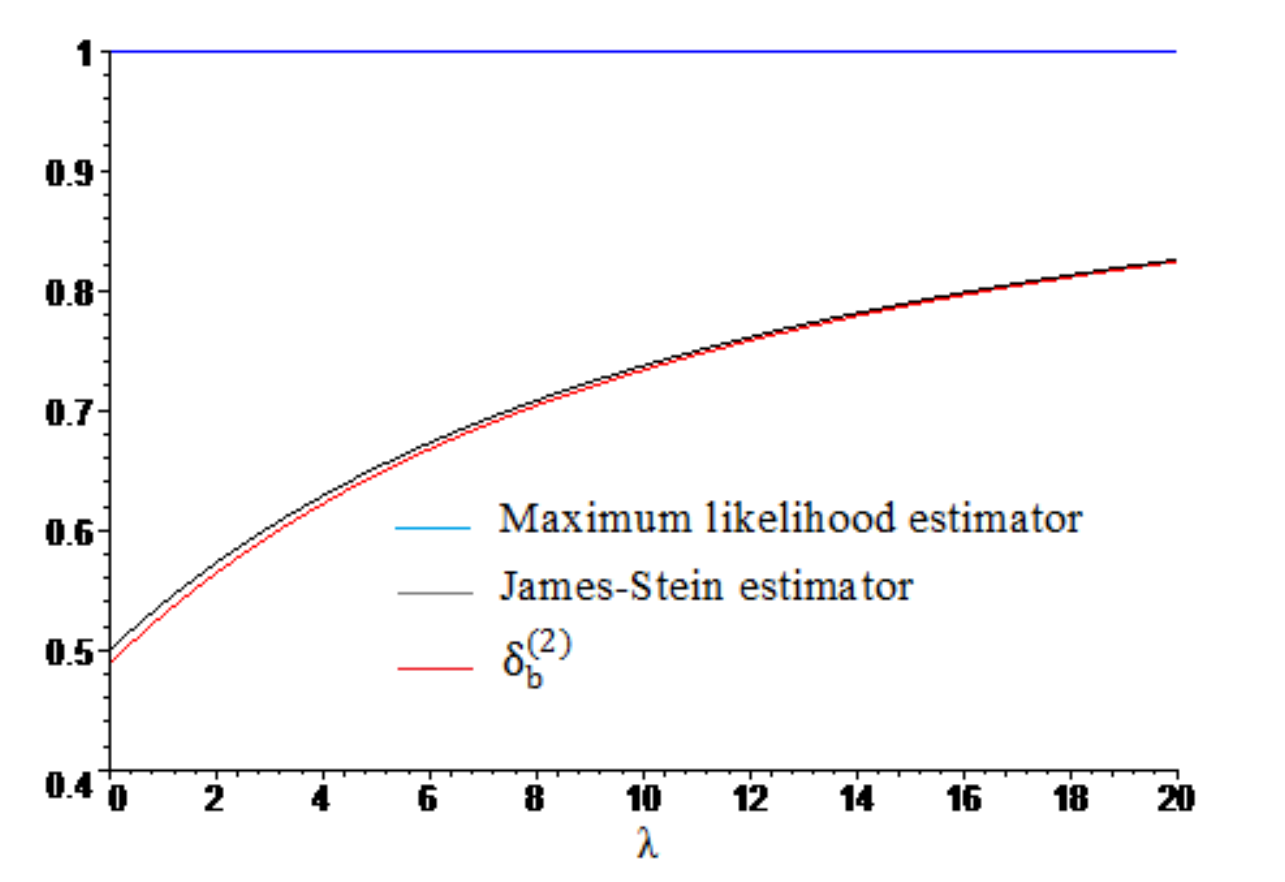}
 \caption{\rm Graph of risks ratios
$\frac{R_{\omega}(\delta_{JS},\theta)}{R_{\omega}(X,\theta)}$ and
$\frac{R_{\omega}(\delta_{b}^{(2)},\theta)}{R_{\omega}(X,\theta)}$ as function of
$\lambda$ for $p=12$ and $\omega=0.4$} \label{figure:f4}
 \centering
 \includegraphics[height=7.1cm,width=10.0cm]{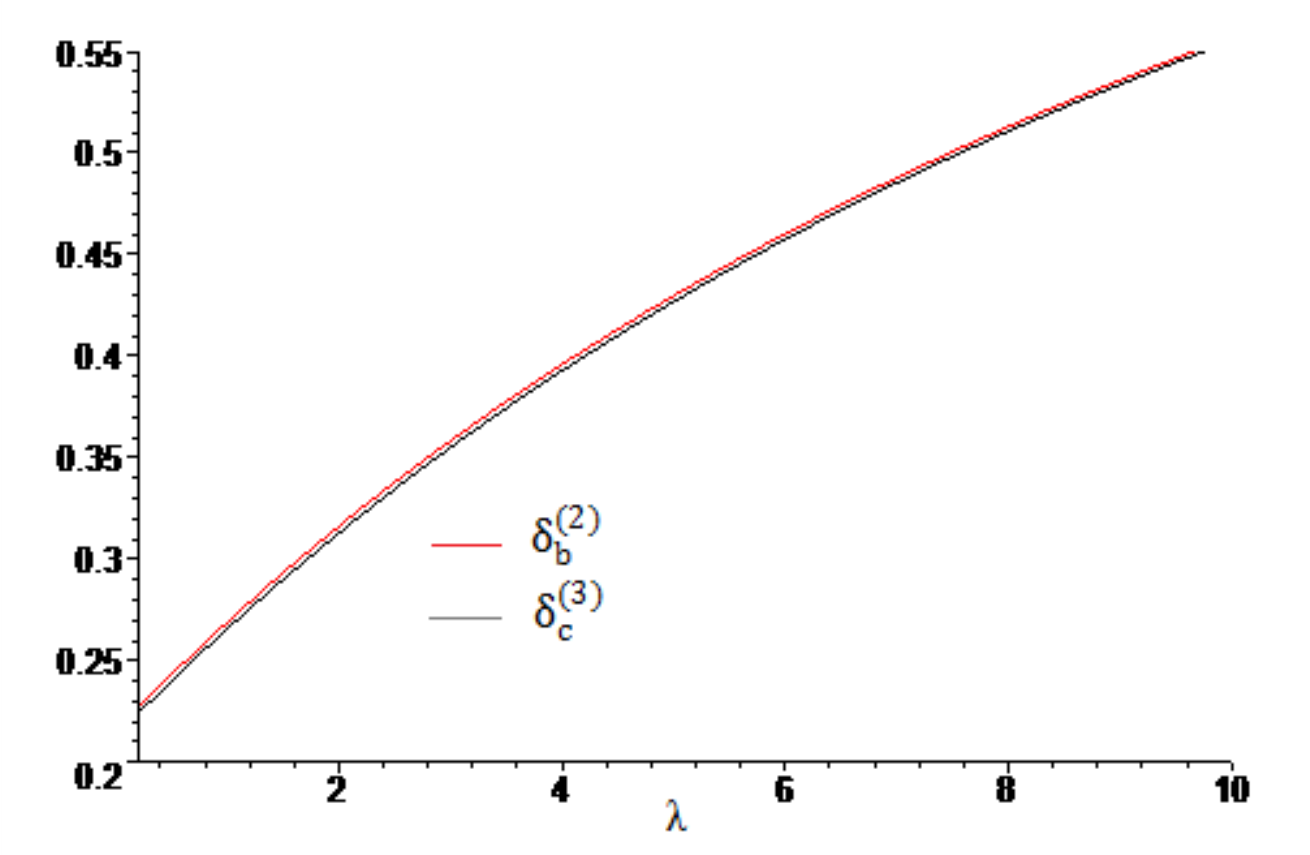}
 \caption{\rm Graph of risks ratios
$\frac{R_{\omega}(\delta_{b}^{(2)},\theta)}{R_{\omega}(X,\theta)}$  and
$\frac{R_{\omega}(\delta_{c}^{(3)},\theta)}{R_{\omega}(X,\theta)}$ as function of
$\lambda$ for $p=14$ and $\omega=0.1$ } \label{figure:f5}
\end{figure}

\newpage

\begin{figure}[htb]
 \centering
 \includegraphics[height=7.1cm,width=10.0cm]{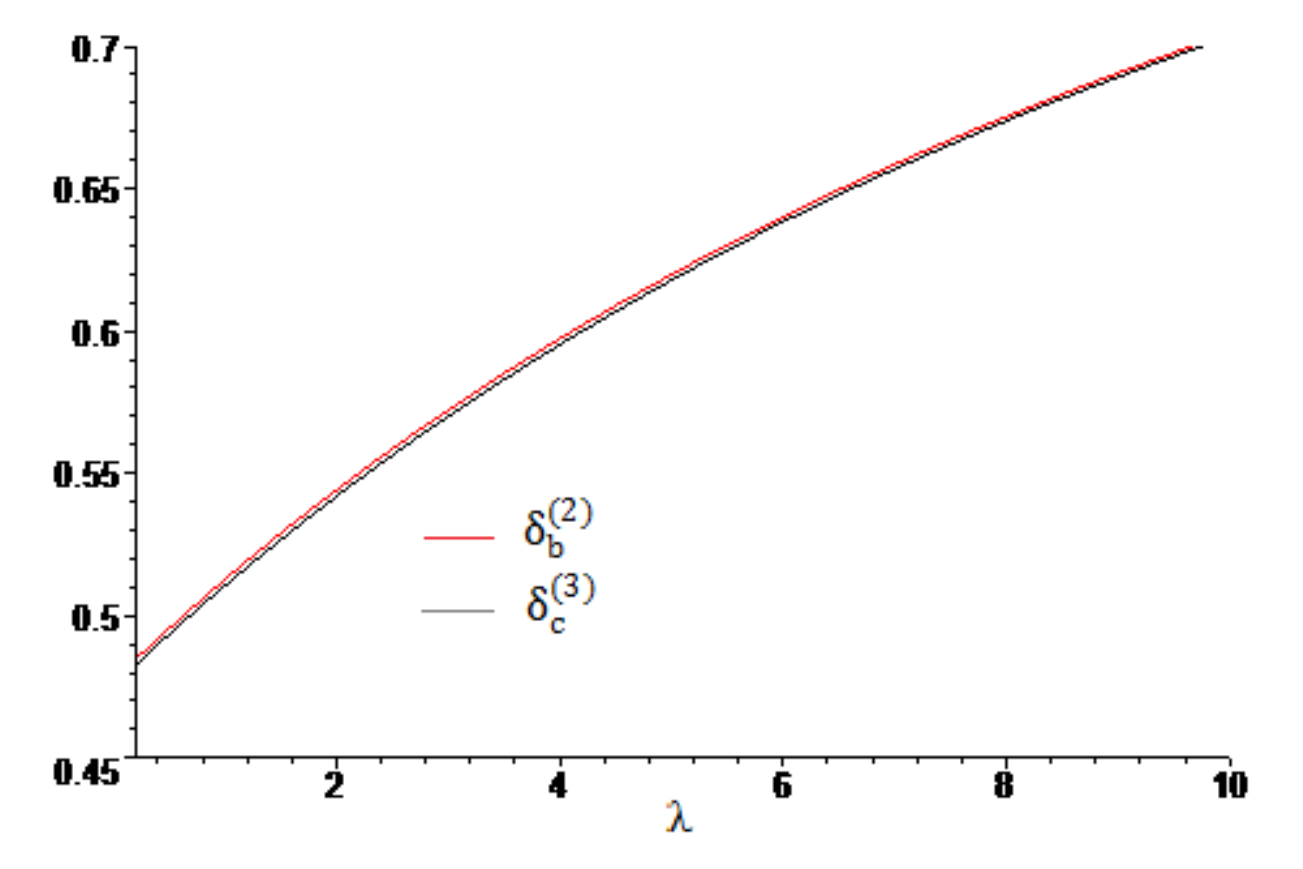}
 \caption{\rm Graph of risks ratios
$\frac{R_{\omega}(\delta_{b}^{(2)},\theta)}{R_{\omega}(X,\theta)}$  and
$\frac{R_{\omega}(\delta_{c}^{(3)},\theta)}{R_{\omega}(X,\theta)}$ as function of
$\lambda$ for $p=14$ and $\omega=0.4$} \label{figure:f6}
 \centering
 \includegraphics[height=7.1cm,width=10.0cm]{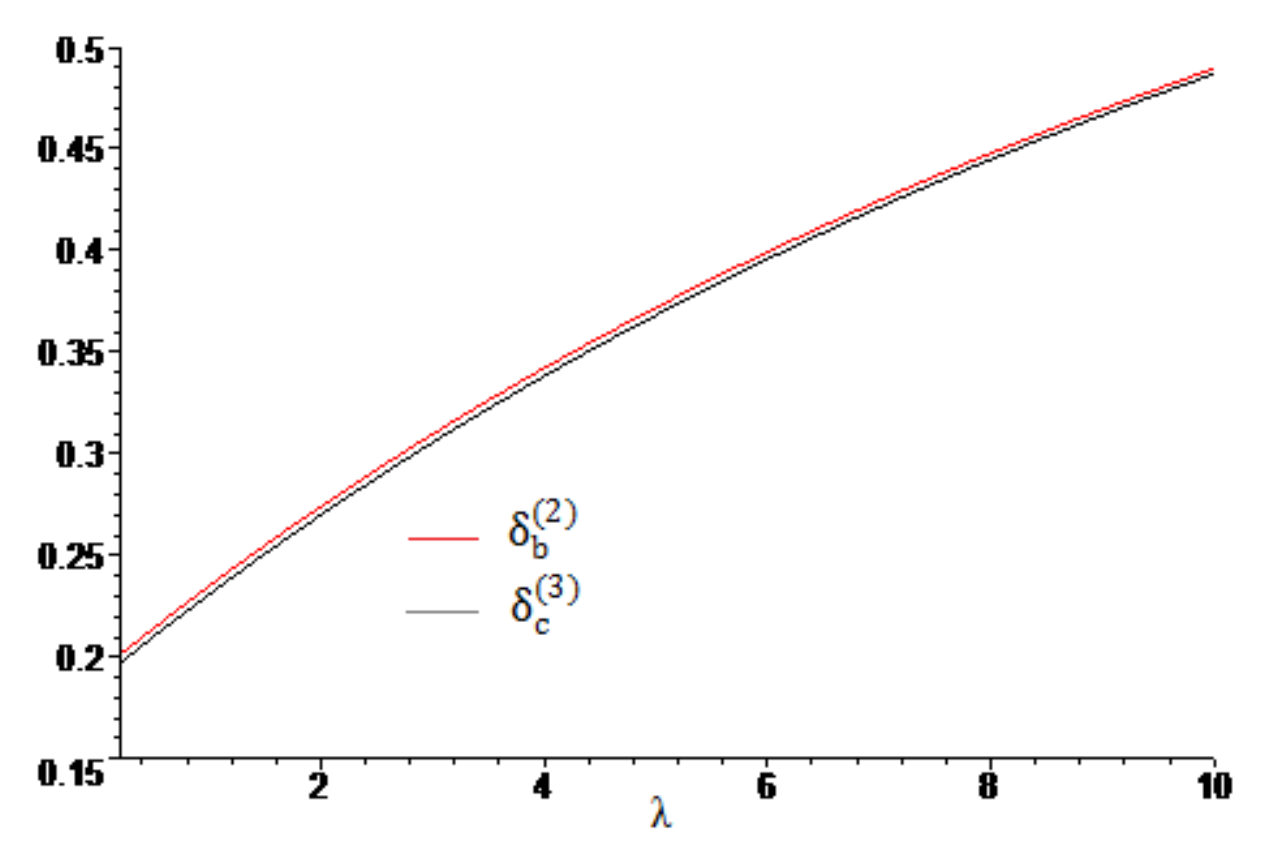}
 \caption{\rm Graph of risks ratios
$\frac{R_{\omega}(\delta_{b}^{(2)},\theta)}{R_{\omega}(X,\theta)}$  and
$\frac{R_{\omega}(\delta_{c}^{(3)},\theta)}{R_{\omega}(X,\theta)}$ as function of
$\lambda$ for $p=18$ and $\omega=0.1$ } \label{figure:f5}
\end{figure}

\newpage

\begin{figure}[htb]
 \centering
 \includegraphics[height=7.1cm,width=10.0cm]{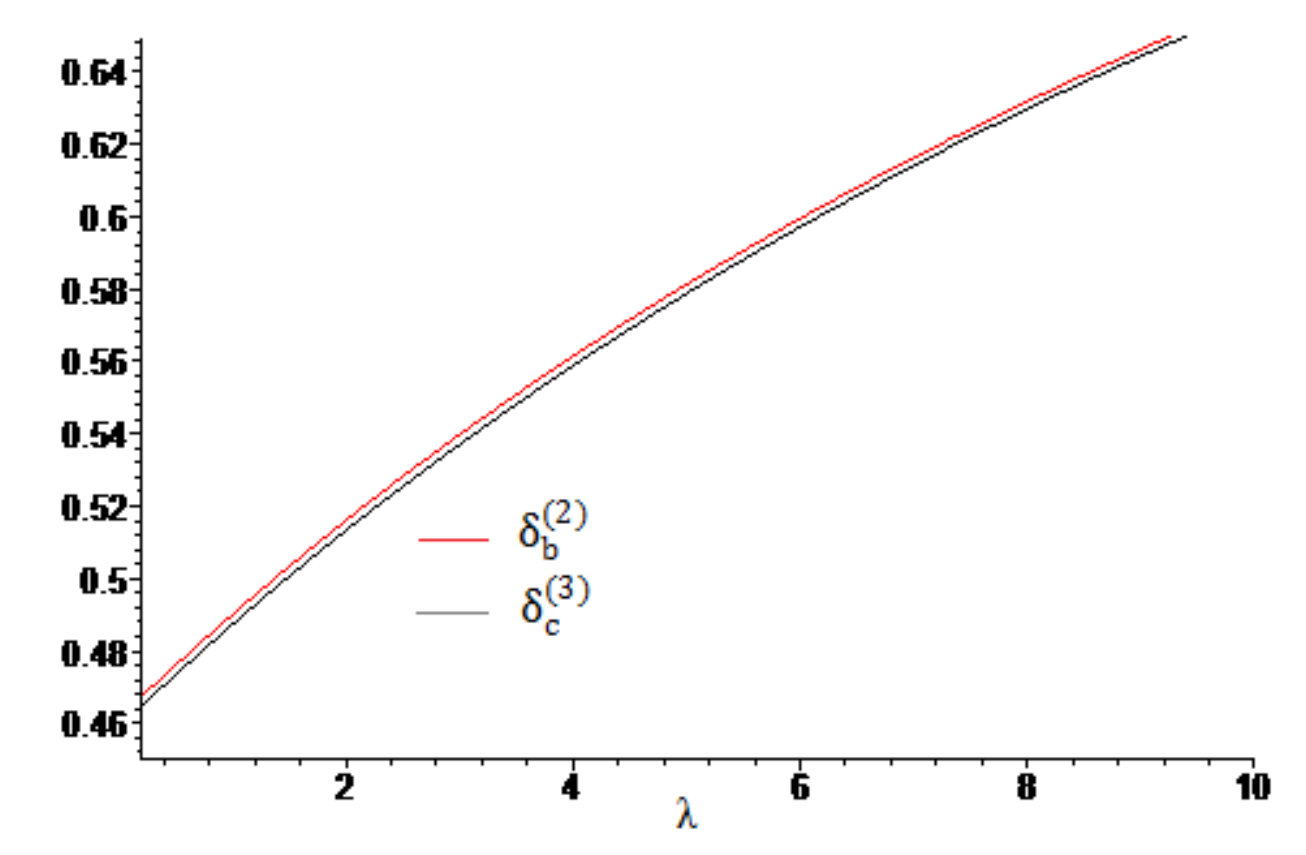}
 \caption{\rm Graph of risks ratios
$\frac{R_{\omega}(\delta_{b}^{(2)},\theta)}{R_{\omega}(X,\theta)}$  and
$\frac{R_{\omega}(\delta_{c}^{(3)},\theta)}{R_{\omega}(X,\theta)}$ as function of
$\lambda$ for $p=18$ and $\omega=0.4$} \label{figure:f6}
\end{figure}

The previous figures show that the risks ratios
$\frac{R_{\omega}(\delta_{JS},\theta)}{R_{\omega}(X,\theta)},$
$\frac{R_{\omega}(\delta_{b}^{(2)},\theta)}{R_{\omega}(X,\theta)}$ and
$\frac{R_{\omega}(\delta_{c}^{(3)},\theta)}{R_{\omega}(X,\theta)}$ are
less than 1, then the estimators $\delta_{JS},\delta_{b}^{(2)},$ and
$\delta_{c}^{(3)}$ dominate the MLE $X$ for divers
values of $p$ and $\omega,$ therefore are minimax. We note that
the estimator $\delta_{b}^{(2)}$ dominates the James-Stein estimator
$\delta_{JS}$ and $\delta_{c}^{(3)}$ dominates $\delta_{b}^{(2)}$ for the selected value of $p$ and $\omega$. We also observe that the gain increases if $\omega$ is near to 0 and decreases if $\omega$ is near to $1.$ The following tables illustrate this note. In these tables, first we give the values of the risks ratios $\frac{R_{\omega}(\delta_{JS},\theta)}{R_{\omega}(X,\theta)},$ $\frac{R_{\omega}(\delta_{b}^{(2)},\theta)}{R_{\omega}(X,\theta)}$ and $\frac{R_{\omega}(\delta_{c}^{(3)},\theta)}{R_{\omega}(X,\theta)}$ for the different values of $\lambda$, $p$ and $\omega$. The first entry is $\frac{R_{\omega}(\delta_{JS},\theta)}{R_{\omega}(X,\theta)},$ the middle entry is $\frac{R_{\omega}(\delta_{b}^{(2)},\theta)}{R_{\omega}(X,\theta)},$ and the third entry is $\frac{R_{\omega}(\delta_{c}^{(3)},\theta)}{R_{\omega}(X,\theta)}$. 

\newpage

\begin{table}[ht]
\caption{$p=14$}
\centering
\begin{tabular}{|C{1.5cm}|C{1.5cm}|C{1.5cm}|C{1.5cm}|C{1.5cm}|C{1.5cm}|C{1.5cm}|}
\hline $\lambda$ & $\omega=0.0$ & $\omega=0.1$ &  $\omega=0.2$ & $\omega=0.5$ & $\omega=0.7$ & $\omega=0.9$ \\
\hline 1.2418 &0.2134  0.2010   0.1973 &0.2920  0.2809  0.2776 &0.3707  0.3608  0.3579 &0.6067 0.6005 0.5987 &0.7640 0.7603 0.7592 &0.9213 0.9201 0.9197 \\
\hline  5.0019 &0.3745  0.3663   0.36309 &0.4371  0.4297  0.4268 &0.4996  0.4930 0.4905 &0.6873 0.6831 0.6815 &0.8124 0.8099 0.8089 &0.9374 0.9366 0.9363 \\
\hline 10.4311 &0.5218  0.5168   0.5150 &0.5697  0.5652  0.5635 &0.6175  0.6135 0.6120 &0.7609 0.7584 0.7575 &0.8565 0.8550 0.8545 &0.9522 0.9517 0.9515 \\
\hline 15.4110 &0.6086  0.6052   0.6041 &0.6477  0.6447  0.6437 &0.6869  0.6842 0.6833 &0.8043 0.8026 0.8020 &0.8826 0.8816 0.8812 &0.9608 0.9605 0.9604 \\
\hline 20.0000 &0.6653  0.6628   0.6621 &0.6988  0.6965  0.6959 &0.7322  0.7302 0.7297 &0.8326 0.8314 0.8310 &0.8996 0.8988 0.8986 &0.9665 0.9663 0.9662 \\
\hline
\end{tabular}
\label{table:nonlin}
\end{table}

\newpage

\begin{table}[ht]
\caption{$p=18$}
\centering
\begin{tabular}{|C{1.5cm}|C{1.5cm}|C{1.5cm}|C{1.5cm}|C{1.5cm}|C{1.5cm}|C{1.5cm}|}
\hline $\lambda$ & $\omega=0.0$ & $\omega=0.1$ &  $\omega=0.2$ & $\omega=0.5$ & $\omega=0.7$ & $\omega=0.9$ \\
\hline  1.2418 &0.1688  0.1608   0.1563 &0.2519  0.2448  0.2406 &0.3351 0.3287 0.3250 &0.5844 0.5804 0.5781 &0.7506 0.7482 0.7469 & 0.9169 0.9161 0.9156 \\
\hline  5.0019 &0.3079  0.3021   0.2980 &0.3771  0.3719  0.3682 &0.4463 0.4417 0.4384 &0.6540 0.6511 0.6490 &0.7924 0.7906 0.7894 & 0.9308 0.9302 0.9298 \\
\hline 10.4311 &0.4535  0.4418   0.4390 &0.5011  0.4976  0.4951 &0.5565 0.5534 0.5512 &0.7228 0.7209 0.7195 &0.8337 0.8325 0.8317 & 0.9446 0.9442 0.9439 \\
\hline 15.4110 &0.5327  0.5299   0.5280 &0.5794  0.5769  0.5752 &0.6261 0.6239 0.6224 &0.7663 0.7649 0.7640 &0.8598 0.8590 0.8584 & 0.9533 0.9530 0.9528 \\
\hline 20.0000 &0.5923  0.5901   0.5888 &0.6331  0.6311  0.6299 &0.6738 0.6721 0.6710 &0.7961 0.7951 0.7944 &0.8777 0.8770 0.8766 & 0.9592 0.9590 0.9589 \\
\hline
\end{tabular}
\label{table:nonlin}
\end{table}

In tables 1 and 2, we note that: if $\omega$ and
$\lambda=\|\theta\|^{2}$ are small, the gain of
the risks ratios
$\frac{R_{\omega}(\delta_{JS},\theta)}{R_{\omega}(X,\theta)},$
$\frac{R_{\omega}(\delta_{b}^{(2)},\theta)}{R_{\omega}(X,\theta)}$
and
$\frac{R_{\omega}(\delta_{c}^{(3)},\theta)}{R_{\omega}(X,\theta)}$
is very important. Also, if the values of $\omega$ and $\lambda$
increase, the gain decreases and approaches to zero, a little
improvement is then obtained. We also observe that, if the values of
$p$ increase, the gain increases and this for each fixed value of
$\omega$. We also see that, if the
values of $p$ are large, the gain is large and consequently
we obtain more improvement. We conclude that, the gain is important
when the parameters $p$ and $\lambda$ are large and $\omega$ is
near to $0$. As seen above, the gain of the risks ratios is
influenced by various values of $p$, $\omega$ and $\lambda$.

Now, we give the tables that present the values of risks
ratios $\frac{R_{\omega}\left(\delta_{
c}^{(3)},\theta\right)}{R_{\omega}(X,\theta)}$ and
$\frac{R_{\omega}\left(\delta_{
d}^{(4)},\theta\right)}{R_{\omega}(X,\theta)}$ for various values of $\lambda$, $p$ and $\omega.$ The first entry is
$\frac{R_{\omega}(\delta_{
c}^{(3)},\theta)}{R_{\omega}(X,\theta)},$ and the second entry is
$\frac{R_{\omega}(\delta_{
d}^{(4)},\theta)}{R_{\omega}(X,\theta)}$.
\newpage
\begin{table}[ht]
\caption{$ p=20$}
\centering
\begin{tabular}{|C{1.5cm}|C{1.5cm}|C{1.5cm}|C{1.5cm}|C{1.5cm}|C{1.5cm}|C{1.5cm}|}
\hline $\lambda$ & $\omega=0.0$ & $\omega=0.1$ &  $\omega=0.2$ & $\omega=0.5$ & $\omega=0.7$ & $\omega=0.9$ \\
\hline 1.2418 &0.1419  0.1414   &0.2277  0.2274  &0.3135 0.3134 &0.5709 0.5713 &0.7426 0.7432 &0.9142 0.9152  \\
\hline 5.0019 &0.2738  0.2732   &0.3464  0.3459  &0.4190 0.4187 &0.6369 0.6368 &0.7821 0.7822 &0.9274 0.9277  \\
\hline10.4311 &0.4091  0.4087   &0.4682  0.4679  &0.5273 0.5270 &0.7045 0.7044 &0.8227 0.8227 &0.9409 0.9410  \\
\hline15.4110 &0.4969  0.4967   &0.5472  0.5470  &0.5975 0.5973 &0.7484 0.7483 &0.8491 0.8490 &0.9497 0.9497  \\
\hline20.0000 &0.5581  0.5579   &0.6022  0.6021  &0.6464 0.6463 &0.7790 0.7790 &0.8674 0.8674 &0.9558 0.9558  \\
\hline
\end{tabular}
\label{table:nonlin}
\end{table}
\begin{table}[ht]
\caption{$ p=24$}
\centering
\begin{tabular}{|C{1.5cm}|C{1.5cm}|C{1.5cm}|C{1.5cm}|C{1.5cm}|C{1.5cm}|C{1.5cm}|}
\hline $\lambda$ & $\omega=0.0$ & $\omega=0.1$ &  $\omega=0.2$ & $\omega=0.5$ & $\omega=0.7$ & $\omega=0.9$ \\
\hline 1.2418 &0.1201  0.1191   &0.2081  0.2074  &0.2961 0.2957 &0.5600 0.5606 &0.7360 0.7372 &0.9120 0.9138  \\
\hline 5.0019 &0.2359  0.2348   &0.3123  0.3114  &0.3887 0.3880 &0.6180 0.6178 &0.7708 0.7710 &0.9236 0.9242  \\
\hline10.4311 &0.3604  0.3596   &0.4244  0.4237  &0.4883 0.4877 &0.6802 0.6799 &0.8081 0.8081 &0.9360 0.9362  \\
\hline15.4110 &0.4448  0.4442   &0.5003  0.4998  &0.5558 0.5554 &0.7224 0.7222 &0.8334 0.8333 &0.9445 0.9445  \\
\hline20.0000 &0.5055  0.5051   &0.5549  0.5546  &0.6044 0.6041 &0.7527 0.7526 &0.8516 0.8516 &0.9505 0.9505  \\
\hline
\end{tabular}
\label{table:nonlin}
\end{table}
\newpage
In tables 3 and 4, we note that: the gain are less than the gain in the
tables 1 and 2, namely there is a little improvement in the domination
of the estimator $\delta_{d}^{(4)}$ to the estimator
$\delta_{c}^{(3)}$ if comparing with the improvement of the
estimator $\delta_{b}^{(2)}$ to the James-Stein estimator or the
improvement of the estimator $\delta_{c}^{(3)}$ to the estimator
$\delta_{b}^{(2)}$. We can also remark that the parameters $p$, $\omega$
and $\lambda$ have the same influence to the risks ratios, as in the tables 1 and 2.
\section{Appendix}
\paragraph{Proof} (Proof of Lemma \ref{l 3.1}) First, we show that, for any real
$\upsilon$
\begin{eqnarray*} \label{e 6.1}
\frac{\partial}{\partial\lambda}E(U^{\upsilon})=\frac{\partial}{\partial\lambda}\int_{R_{+}}x^{\upsilon}\chi_{p}^{2}(\lambda;dx)
=\upsilon 2^{\upsilon-1}\sum^{+\infty}_{k=0}
\frac{\Gamma(\frac{p}{2}+\upsilon+k)}{\Gamma(\frac{p}{2}+1+k)}P\left(\frac{\lambda}{2};dk\right),
\end{eqnarray*}
where $P(\frac{\lambda}{2})$ being the Poisson distribution of
parameter $\frac{\lambda}{2}$.\\
Using the formula (\ref{e 2}) we have, for any real $\upsilon$
\begin{eqnarray} \label{e 6.2}
E(U^{\upsilon})=E[(\chi_{p}^{2}(\lambda))^{\upsilon}]=E[(\chi_{p+2K}^{2})^{\upsilon}]
=2^{\upsilon}E\left[\frac{\Gamma(\frac{p}{2}+K+\upsilon)}{\Gamma(\frac{p}{2}+K)}\right],
\end{eqnarray}
where $K \sim P(\frac{\lambda}{2})$. Then
\begin{eqnarray*}
\frac{\partial}{\partial\lambda}E(U^{\upsilon})&=&\frac{\partial}{\partial\lambda}\int_{R_{+}}x^{\upsilon}\chi_{p}^{2}
(\lambda;dx)\\
&=&2^{\upsilon}\sum^{+\infty}_{k=0}\left[\frac{\Gamma(\frac{p}{2}+k+\upsilon)}{\Gamma(\frac{p}{2}+k)}\right]\frac{1}{k!}
\frac{\partial}{\partial\lambda}\left[\left(\frac{\lambda}{2}\right)^{k}exp\left(-\frac{\lambda}{2}\right)\right]\\
&=&2^{\upsilon-1}\sum^{+\infty}_{k=0}\left[\frac{\Gamma(\frac{p}{2}+k+\upsilon)}{\Gamma(\frac{p}{2}+k)}\right]\frac{1}{k!}
exp\left(-\frac{\lambda}{2}\right)\left[-\left(\frac{\lambda}{2}\right)^{k}+k\left(\frac{\lambda}{2}\right)^{k-1}\right]\\
&=&2^{\upsilon-1}exp\left(-\frac{\lambda}{2}\right)\left\{-\sum^{+\infty}_{k=0}\left[\frac{\Gamma(\frac{p}{2}+k+\upsilon)}
{\Gamma(\frac{p}{2}+k)}\right]\frac{1}{k!}\left(\frac{\lambda}{2}\right)^{k}\right\}\\
&+&2^{\upsilon-1}exp\left(-\frac{\lambda}{2}\right)\left\{\sum^{+\infty}_{k=0}\left[\frac{\Gamma(\frac{p}{2}+k+\upsilon+1)}{\Gamma(\frac{p}{2}+k+1)}\right]\frac{1}
{k!}\left(\frac{\lambda}{2}\right)^{k}\right\}\\
\end{eqnarray*}
\begin{eqnarray*}
&=&2^{\upsilon-1}exp\left(-\frac{\lambda}{2}\right)\left\{\sum^{+\infty}_{k=0}\frac{1}
{k!}\left(\frac{\lambda}{2}\right)^{k}\left[\frac{\Gamma(\frac{p}{2}+k+\upsilon)}{\Gamma(\frac{p}{2}+k+1)}\right]\left[-\left(\frac{p}{2}+k\right)+\left(\frac{p}{2}+\upsilon+k\right)\right]\right\}\\
&=&\upsilon 2^{\upsilon-1}\sum^{+\infty}_{k=0}
\frac{\Gamma(\frac{p}{2}+\upsilon+k)}{\Gamma(\frac{p}{2}+1+k)}P\left(\frac{\lambda}{2};dk\right).
\end{eqnarray*}
Let the function
\begin{eqnarray*}
K_{p,r,s}(\lambda)&=&\left(\frac{\partial}{\partial\lambda}
\int_{R_{+}}x^{r}\chi_{p}^{2}(\lambda;dx)\right)\left(
\int_{R_{+}}x^{s}\chi_{p}^{2}(\lambda;dx)\right)\\
&-&\left(\frac{\partial}{\partial\lambda}
\int_{R_{+}}x^{s}\chi_{p}^{2}(\lambda;dx)\right)\left(\int_{R_{+}}x^{r}\chi_{p}^{2}(\lambda;dx)\right).
\end{eqnarray*}
For the function $H_{p,r,s}$ to be strictly increasing, it suffices
that the function $K_{p,r,s}$ takes positive values. From the
equality (\ref{e 6.2}), we obtain
\begin{eqnarray*}
K_{p,r,s}(\lambda)&=&2^{r+s-1}r\sum^{+\infty}_{i=0}\sum^{+\infty}_{j=0}
\frac{\Gamma(\frac{p}{2}+r+i)}{\Gamma(\frac{p}{2}+i+1)}\frac{\Gamma(\frac{p}{2}+s+j)}
{\Gamma(\frac{p}{2}+j)}P\left(\frac{\lambda}{2};di\right)P\left(\frac{\lambda}{2};dj\right)\\
&-&2^{r+s-1}s\sum^{+\infty}_{i=0}\sum^{+\infty}_{j=0}
\frac{\Gamma(\frac{p}{2}+r+j)}{\Gamma(\frac{p}{2}+j)}\frac{\Gamma(\frac{p}{2}+s+i)}
{\Gamma(\frac{p}{2}+i+1)}P\left(\frac{\lambda}{2};dj\right)P\left(\frac{\lambda}{2};di\right).
\end{eqnarray*}
As, $r>s$ then
\begin{eqnarray*}
K_{p,r,s}(\lambda)&\geq& r2^{r+s-1}\sum^{+\infty}_{i=0}\sum^{+\infty}_{j=0}
l_{p,r,s}(i,j)P\left(\frac{\lambda}{2};di\right)P\left(\frac{\lambda}{2};dj\right),
\end{eqnarray*}
where
\begin{eqnarray*}
l_{p,r,s}(i,j)=\frac{\Gamma(\frac{p}{2}+r+i)\Gamma(\frac{p}{2}+s+j)-\Gamma(\frac{p}{2}+r+j)\Gamma(\frac{p}{2}+s+i)}
{\Gamma(\frac{p}{2}+i+1)\Gamma(\frac{p}{2}+j)}.
\end{eqnarray*}
We note that, for any $i,$  $l_{p,r,s}(i,j)=0;$ then we have
\begin{eqnarray*}
K_{p,r,s}(i,j)&\geq& r2^{r+s-1}\sum^{+\infty}_{i=0}\sum^{+\infty}_{j>i}
(l_{p,r,s}(i,j)+l_{p,r,s}(j,i))P\left(\frac{\lambda}{2};di\right)P\left(\frac{\lambda}{2};dj\right).
\end{eqnarray*}
But if $i<j$, we get
\begin{eqnarray*}
l_{p,r,s}(i,j)+l_{p,r,s}(j,i)&=&\left(\Gamma\left(\frac{p}{2}+r+i\right)\Gamma\left(\frac{p}{2}+s+j\right)
-\Gamma\left(\frac{p}{2}+r+j\right)\Gamma\left(\frac{p}{2}+s+i\right)\right)\\
& \times&\left[\frac{1}{\Gamma(\frac{p}{2}+i+1)\Gamma(\frac{p}{2}+j)}-\frac{1}{\Gamma(\frac{p}{2}+j+1)\Gamma(\frac{p}{2}+i)}\right]\\
&=&\frac{\Gamma(\frac{p}{2}+r+i)\Gamma(\frac{p}{2}+s+i)}{\Gamma(\frac{p}{2}+i)\Gamma(\frac{p}{2}+j)}
\left[\frac{1}{\frac{p}{2}+i}-\frac{1}{\frac{p}{2}+j}\right]\\
&\times &\left[\prod^{j-i-1}_{t=0}\left(\frac{p}{2}+s+i+t\right)-\prod^{j+i-1}_{t=0}\left(\frac{p}{2}+r+i+t\right)\right]\\
&\leq&0,
\end{eqnarray*}
because for any $t,$ $\frac{p}{2}+s+i+t<\frac{p}{2}+r+i+t$. As in hypothesis $r<0,$ we
have $K_{p,r,s}(\lambda)>0.$ Thus, we obtain the desired result.\\
ii) Using i) it is clear that the function
$H_{p,r}^{1}(\lambda)=\frac{E(\|X\|^{-r})}{E(\|X\|^{-2r+2})}$ is non-decreasing on $\lambda$, then the function $\frac{1}{H_{p,r}^{1}(\lambda)}$ is non-increasing on $\lambda$, thus
\begin{eqnarray*}
\sup_{\|\theta\|}\left(\frac{E(\|X\|^{-2r+2})}{E(\|X\|^{-r})}\right)&=&\sup_{\|\theta\|}\left(\frac{1}{H_{p,r}^{1}(\lambda)}\right)\\
&=&\frac{1}{H_{p,r}^{1}(0)}\\
&=&2^{\frac{-r+2}{2}} \frac{\Gamma(\frac{p}{2}-r+1)}{\Gamma(\frac{p-r}{2})}.
\end{eqnarray*}\cqfd\\
\section*{Conclusion}
In this work, we studied the estimating of the the mean $\theta$ of
a multivariate normal distribution $X\sim N_{p}\left(\theta,
\sigma^{2}I_{p}\right)$ where $\sigma^{2}$ is known. The criterion
adopted for comparing two estimators is the risk associated to the
balanced loss function. First, we established the minimaxity of the
estimators defined by
$\delta_{a}^{(1)}=\left(1-a/\|X\|^{2})\right)X$, where the real
parameter $a$ may depend on $p$ and we constructed the
James-Stein estimator that has the minimal risk in this class.
Secondly, we considered the estimators of polynomial form with the
indeterminate $1/\|X\|^{2}$ and showed that if we increase
the degree of the polynomial, we can construct a better estimator.
We concluded that we constructed a series of estimators of
polynomial form such that if we increase the degree, the estimator
becomes much better. An extension of this work is to obtain the
similar results in the case where the model has a symmetrical
spherical distribution.

\end{document}